\documentclass[11pt]{article}
\usepackage{amsmath,amsfonts,latexsym,amssymb,chicago,graphicx,annals,here}
\advance \topmargin by -\headheight \advance \topmargin by
-\headsep

 \input epsf 

\newcommand{\mycite}[1]{{\small \sc \citeNP{#1}}}

\def\P{{\mathbb P}}

\def\R{\mathbb R}

\def\ai{\mbox{Ai}}
\def\bi{\mbox{Bi}}
\def\a{\alpha}

\def\eop{\hfill\mbox{$\Box$}\newline}
\newtheorem{theorem}{Theorem}[section]
\newtheorem{corollary}{Corollary}[section]
\newtheorem{lemma}{Lemma}[section]

\newtheorem{remark}{Remark}[section]

\def\bigO{{\cal O}}
\def\calP{{\cal P}}
\def\phase{{\rm ph}}

\begin{document}
\title{The tail of the maximum of Brownian motion minus a parabola}
\author{Piet Groeneboom and Nico M.\ Temme}
\date{\today}
\affiliation{Delft University of Technology, Mekelweg 4, 2628 CD Delft, The Netherlands,
p.groeneboom@tudelft.nl; http://dutiosc.twi.tudelft.nl/$_{\widetilde{~}}$pietg/\\
and CWI, Science Park 123, 1098 XG Amsterdam, The Netherlands,\\
Nico.Temme@cwi.nl; http://homepages.cwi.nl/$_{\widetilde{~}}$nicot/}
\AMSsubject{60J65,60J75}
\maketitle
\begin{abstract}
We analyze the tail behavior of the maximum $N$ of $\{W(t)-t^2:t\ge0\}$, where $W$ is standard Brownian motion on $[0,\infty)$ and give an asymptotic expansion for $\P\{N\ge x\}$, as $x\to\infty$. This extends a first order result on the tail behavior, which can be deduced from \mycite{huesler:99}. We also point out the relation between certain results in \mycite{piet:10} and \mycite{janson:10}.
\end{abstract}

\section{Introduction}
\label{section:intro}
The distribution function of the maximum of Brownian motion minus a parabola was studied in the two recent papers \mycite{janson:10} and \mycite{piet:10}, both for one-sided and two-sided Brownian motion. The characterization of the distribution function is somewhat different in the two papers, but both characterizations (unavoidably) involve Airy functions. In this note we address the tail behavior of the distribution, a topic that was not addressed in these papers.

The tail behavior of the maximum plays an important role in certain recent studies on the asymptotic distribution of tests for monotone hazards, based on integral-type statistics measuring the distance between the empirical cumulative hazard function and its greatest convex minorant, for example in \mycite{GrJo:10}.

Let $N$ be defined by
\begin{equation}
\label{def_max}
N=\max_{t\ge0}\{W(t)-t^2\},\,t\ge0,
\end{equation}
where $W$ is standard Brownian motion on $[0,\infty)$. It can be deduced from Theorem 2.1 in\\ \mycite{huesler:99} that the distribution function $F_N$ of $N$ satisfies:
$$
1-F_N(x)\sim \frac1{\sqrt{3}}\exp\left\{-\frac{8x^{3/2}}{3\sqrt{3}}\right\},\,x\to\infty.
$$
In section 2 we will give an asymptotic expansion of the left-hand side, which extends this result.
The proof is based on an integral expression for the density, derived from \mycite{piet:10} (which in turn relies on \mycite{piet:89}), and uses a saddle point method for the integral over a shifted path in the complex plane. As a side effect, it also leads to a clarification of the relation between the representations of the distribution, given in \mycite{janson:10} and \mycite{piet:10}.

\section{Main results}
\label{section:main_result}
In the following, we will use Corollary 2.1 of \mycite{piet:10}, which is stated below for ease of reference, specialized to the density of the maximum of $W(t)-t^2$ (instead of the more general $W(t)-ct^2$).

\begin{lemma}
{\rm (Corollary 2.1 in \mycite{piet:10})}
\label{lemma:dens_M}
The density $f$ of $N$ is given by:
\begin{equation}
\label{density_representation}
f_N(x)=2^{2/3}\left\{{\mathrm {Ai}}\bigl(2^{2/3}x\bigr) -2\,{\mathrm {Re}}\left(e^{-i\pi/6}\int_0^{\infty}\frac{{\rm {\mathrm {Ai}}}\left(e^{-i\pi/6}u\right){\mathrm {Ai}}\,'\bigl(iu+2^{2/3}x\bigr)}{\mathrm {Ai}(iu)}\,du\right)\right\},\,x>0.
\end{equation}
where ${\mathrm {Ai}}$ is the Airy function ${\mathrm {Ai}}$, as defined in, e.g., \mycite{olver:10}\footnote{http://dlmf.nist.gov/9}.
\end{lemma}

We deduce from this the following representation which is better suited for our purposes.

\begin{lemma}
\label{lemma:dens_representation}
The density $f_N$ of $N$ is given by:
\begin{equation}
\label{density_representation2}
f_N(x)=\frac{2^{2/3}}{\pi}\,{{\mathrm {Re}}}\left(\int_0^{\infty}\frac{\mathrm {Ai}\bigl(iu+2^{2/3}x\bigr)}{\mathrm {Ai}(iu)^2}\,du\right)=\frac1{2^{1/3}\pi}\int_{-\infty}^{\infty}\frac{\mathrm {Ai}\bigl(iu+2^{2/3}x\bigr)}{\mathrm {Ai}(iu)^2}\,du,\,x>0.
\end{equation}
\end{lemma}

\noindent
{\bf Proof.} Integration by parts of the second term of (\ref{density_representation}) yields:
$$
f_N(x)=2\cdot2^{2/3}\,\mbox{\rm Re}\left(\int_0^{\infty}\ai\bigl(iu+2^{2/3}x\bigr)\frac{d}{du}\left\{\frac{e^{-i\pi/6}\ai\left(e^{-i\pi/6}u\right)}{i\ai(iu)}\right\}\,du\right),\,x>0.
$$
Let the function $h$ be defined by
$$
h(u)=\frac{d}{du}\left\{\frac{e^{-i\pi/6}\ai\left(e^{-i\pi/6}u\right)}{i\ai(iu)}\right\}.
$$
Using  \mycite{olver:10}\footnote{http://dlmf.nist.gov/9.2.E11}:
$$
\ai\left(e^{-i\pi/6}u\right)=\ai\left(e^{-2i\pi/3}iu\right)=\tfrac12e^{-i\pi/3}\left\{\ai(iu)+i\bi(iu)\right\},
$$
we obtain
$$
h(u)=\tfrac12\left\{\frac{\ai\,'(iu)+i\bi\,'(iu)}{i\ai(iu)}
-\frac{\left(\ai(iu)+i\bi(iu)\right)\ai\,'(iu)}{i\ai(iu)^2}\right\},
$$
and using the Wronskian $\ai(z)\bi\,'(z)-\ai\,'(z)\bi(z)=1/\pi$ we conclude that
$$
h(u)=\frac{1}{2\pi\ai(iu)^2}\,.
$$
This gives the desired result.
\eop

\begin{remark}
{\rm Lemma \ref{lemma:dens_representation} is in fact equivalent to relation (5.10) in \mycite{janson:10}. The difference in the scaling constants is caused by the fact that they consider the maximum of $W(t)-\tfrac12t^2$ instead of the maximum of $W(t)-t^2$ (see also section \ref{section:conclusion}) and the fact that they integrate from $-\infty$ to $\infty$ (in that way also the imaginary part drops out). However, they arrive at this relation in a completely different way. So in this case we can go from Corollary 2.1 in \mycite{piet:10} to the result in \mycite{janson:10}, just by using integration by parts. This might serve as a first step in establishing the relation between the representations in the two papers.
}
\end{remark}

We are now ready to prove our main result. We will give two proofs, one based on the first equality in (\ref{density_representation2}) and the other one based on the second equality.

\begin{theorem}
\label{th:main_result}
Let $N$ be defined by (\ref{def_max}), and let $f_N$ and $F_N$ be the density and the distribution function of $N$, respectively. Then:
\begin{enumerate}
\item[(i)]
\begin{equation}
\label{asymp_representation2}
f_N(x)\sim2^{2/3}\sqrt{\pi}\sum_{k=0}^\infty a_k \Phi_k(x),
\end{equation}
where
$$
\Phi_k(x)=\frac1{2\pi i}
\int_{\infty e^{-\pi i/3}}^{\infty e^{\pi i/3}} e^{\frac14t^3-2^{2/3}xt}t^{\frac32-3k}\,dt.
$$
and $a_k$ is given by:
\begin{align*}
\label{trans4}
&a_0=  1,\quad 
a_1=  2,\quad 
a_2=  -20,\quad 
a_3=  560,\quad
a_4=  -25520,\\
&a_5=  1601600,\quad 
a_6=  -127568000,\quad 
a_7=  12287436800.
\end{align*}


\item[(ii)] More explicitly:
\begin{equation}
\label{fNxasymp1}
f_N(x)\sim\frac{4\sqrt{x}}{3}\exp\left(-\frac{8x^{3/2}}{3\sqrt{3}}\right)\sum_{k=0}^\infty \frac{b_k}{x^{\frac32k}},\quad x\to\infty,
\end{equation}
where the first coefficients are
$$
b_0=  1,\quad 
b_1=  \tfrac{19}{48}\sqrt{3},
\quad 
b_2=  -\tfrac{3851}{1536},\quad 
b_3=   \tfrac{3380005}{221184}\sqrt{3},\quad
b_4=-\tfrac{6474441455}{14155776}.
$$
\item[(iii)]
$$
1-F_N(x)\sim \frac1{\sqrt{3}}\exp\left\{-\frac{8x^{3/2}}{3\sqrt{3}}\right\}\sum_{k=0}^{\infty}\frac{c_k}{x^{\tfrac32k}},\,x\to\infty.,
$$
where the first coefficients are:
$$
c_0= 1, \quad
c_1= \tfrac{19}{48}\sqrt{3},\quad
c_2= -\tfrac{4535}{1536},\quad
c_3= \tfrac{3869785}{221184}\sqrt{3},\quad
c_4= -\tfrac{7310315015}{14155776}.
$$
\end{enumerate}
\end{theorem}
\noindent
{\bf Proof}.
Here we only derive the leading terms. Further terms in the asymptotic expansion are computed in the appendix.

Let $g(t)$ be defined by
\begin{equation}\label{gtdef}
g(t)=\int_0^{\infty}\frac{e^{-itu}}{\pi\ai(iu)^2}\,du.
\end{equation}
By using the representation
$$
\ai(z)=\frac1{2\pi i}\int_{\infty e^{-\pi i/3}}^{\infty e^{\pi i/3}}
\exp\left\{\tfrac13t^3-zt\right\}\,dt.
$$
the second integral representation of $f_N(x)$ in \eqref{density_representation2} becomes
\begin{equation}
\label{asymp_representation}
f_N(x)=2^{-1/3}\mbox{\rm Re}\left\{\frac1{\pi i}\int_{\infty e^{-\pi i/3}}^{\infty e^{\pi i/3}} \exp\left\{\tfrac13t^3-2^{2/3}xt\right\}g(t)\,dt\right\},
\end{equation}
and this integral will be expanded for large values of $x$. For this we need an expansion of $g(t)$ for large values of $t$.

For large values of the argument $|u|$, and $|\mbox{arg}(iu)|<\pi$, we have (see \mycite{olver:10}\footnote{http://dlmf.nist.gov/9.7} for the asymptotic behavior of the Airy functions):
\begin{equation}
\label{inside_function}
\frac1{\pi\ai(iu)^2}\sim 2\exp\left\{\tfrac43(iu)^{3/2}\right\}
(1+i)\sqrt{2u}.
\end{equation}
Hence, the dominant part of the integral in \eqref{gtdef} is $e^{\phi(u)}$, where
\begin{equation}\label{phidef}
\phi(u)=\tfrac43(iu)^{3/2}-itu.
\end{equation}
The derivative of this function vanishes at $u_s=-it^2/4$, and we have at $u_s$ the expansion
\begin{equation}\label{phiexpf}
\phi(u)=-\tfrac1{12}t^3-\frac{1}{t}\left(u-u_s\right)^2+\bigO\left(\left(u-u_s\right)^3)/t^3\right).
\end{equation}

This suggests to take $u_s$ as a saddle point (for a changed integration road) for the integral in \eqref{gtdef}.

\begin{figure}
\begin{center}
\epsfxsize=12cm \epsfbox{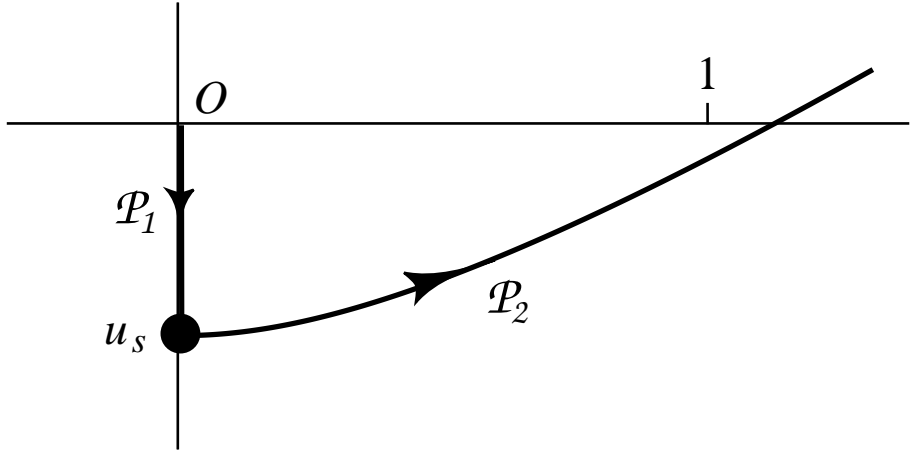}
\caption{Modification of the path of integration for the integral in \eqref{gtdef} from $(0,\infty)$ to $\calP_1\cup\calP_2$, for $t=1$ with saddle point at $u_s=-i/4$. \label{fig1}
}
\end{center}
\end{figure}

We consider the following integration path: first the path $\calP_1$, going from $0$ to $u_s$, and next the path $\calP_2$, from $u_s$ to $+\infty$, into the valley of $\exp(\tfrac43(iu)^{3/2})$. That is, at $+\infty$ the phase of $u$ is $\frac16\pi$. See Figure~\ref{fig1}, where we have shown the paths $\calP_1$ and $\calP_2$  for $t=1$.

We write $g(t)=g_1(t)+g_2(t)$, where $g_j(t)$ is the contribution from the path $\calP_j$, $j=1,2$. Then:
$$
g_1(t)=\int_{\calP_1}\frac{e^{-itu}}{\pi\ai(iu)^2}\,du
=-i\int_0^{t^2/4}\frac{e^{-tu}}{\ai(u)^2}\,du,
$$
This shows that $g_1(t)$ is purely imaginary for $t>0$.  When we replace in \eqref{asymp_representation}  $g(t)$ with $g_1(t)$, we see that this contribution to $f_N(x)$ becomes
$$
-2^{-1/3}\mbox{\rm Re}\left\{\frac1{\pi }\int_{\infty e^{-\pi i/3}}^{\infty  e^{\pi i/3}}G(t)\,dt\right\},
$$
where $G(t)$ is an analytic function which is real for positive values of $t$. The behavior of $G(t)$ at infinity allows to take the contour along the imaginary axis. Integrating in this way, we obtain
$$
-2^{-1/3}\mbox{\rm Re}\left\{\frac{i}{\pi }\int_{0}^{\infty}\left(G(iv)+G(-iv)\right)\,dv\right\},
$$
and since $G(iv)+G(-iv)$ is real for real $v$ (because $G(t)$ is a real function), this integral is purely real, and, hence, we can ignore this contribution.

Next we consider the integral over $\calP_2$. A parametrization of this path follows from the equation 
$\mbox{\rm Im}\{\phi(u)\}= \mbox{\rm Im}\{\phi(u_s)\}$, and we see from \eqref{phiexpf} that  $\mbox{\rm Im}\{\phi(u_s)\}=0$.  By writing $u=r\exp(i\theta)$, it follows that the path $\calP_2$ can be described in terms of polar coordinates by 
$$
r= \frac{3t^2\cos^2\theta}{16\sin^2\left(\frac32\theta+\frac34\pi\right)},\quad -\tfrac12\pi\le \theta\le \tfrac16\pi,
$$
where for $\theta\to-\frac12\pi$ we have to apply l'H\^opital's rule to obtain $r=\frac14t^2$.

There is no need to follow this path for obtaining the asymptotic expansion of $g_2(t)$ and for the sake of convenience we use the path from $u_s$ to $u_s+\infty$ parallel to the positive real axis.  In this way we find the contribution from the path $\calP_2$, that is
\begin{equation}\label{defg2}
g_2(t)=\int_{\calP_2}\frac{e^{-itu}}{\pi\ai(iu)^2}\,du=\int_{u_s} ^{u_s+\infty}\frac{e^{-itu}}{\pi\ai(iu)^2}\,du.
\end{equation}

At the lower limit $u_s=-it^2/4$ of the integration, we have by (\ref{inside_function}),
$$
\frac{\exp\left\{-\tfrac43(iu_s)^{3/2}\right\}}{\pi\ai(iu_s)^2}\sim 2(1+i)\sqrt{2u_s}=2t,
$$
giving
$$
g_2(t)\sim2t\int_{u_s} ^{u_s+\infty}e^{\phi(u)}\,du.
$$
Next we neglect the $\bigO-$term in the expansion in \eqref{phiexpf}, and substitute
$u=u_s+t^2w,\,w\in\R$. This yields
\begin{equation}
\label{Regtasymp}
\mbox{Re}(g(t))=\mbox{Re}(g_2(t))\sim 2t^3e^{-\frac1{12}t^3}\int_0^{\infty}e^{-t^3w^2}\,dw=t^{3/2}e^{-\frac1{12}t^3}\sqrt{\pi}, \quad t\to\infty.
\end{equation}
\eop


Plugging the result in \eqref{Regtasymp} into (\ref{asymp_representation}), 
and expanding at $t=2^{4/3}\sqrt{x/3}$ yields:
\begin{align*}
f_N(x)&\sim \mbox{\rm Re}\left\{\frac{2^{-1/3}}{i\sqrt{\pi}}\int_{\infty e^{-\pi i/3}}^{\infty e^{\pi i/3}} t^{3/2}\exp\left\{\tfrac13t^3-4^{1/3}xt-\tfrac1{12}t^3\right\}\,dt\right\}\\
&=\mbox{\rm Re}\left\{\frac{2^{-1/3}}{i\sqrt{\pi}}\int_{\infty e^{-\pi i/3}}^{\infty e^{\pi i/3}} t^{3/2}\exp\left\{\tfrac14t^3-4^{1/3}xt\right\}\,dt\right\}\\
&\sim\frac{2^{5/3}(x/3)^{3/4}}{\sqrt{\pi}}\exp\left\{-\frac{8x^{3/2}}{3\sqrt{3}}\right\}\int_{-\infty}^{\infty}\exp\left\{-2^{-2/3}u^2\sqrt{3x}\right\}\,du\\
&=\frac{4\sqrt{x}}{3}\exp\left\{-\frac{8x^{3/2}}{3\sqrt{3}}\right\},\quad x\to\infty.
\end{align*}
This gives the leading term of the expansion; the further coefficients are computed in the appendix.
Part (iii) follows by integrating this relation.\\
(Second proof.) Again we only derive the leading term. We start with the second representation in (\ref{density_representation2}).
\begin{equation}
\label{density_representation3}
f_N(x)=\frac1{2^{1/3}\pi}\int_{-\infty}^{\infty}\frac{\mathrm {Ai}\bigl(iu+2^{2/3}x\bigr)}{\mathrm {Ai}(iu)^2}\,du
=\frac1{2^{1/3}\pi i}\int_{-i\infty}^{i\infty}\frac{\mathrm {Ai}\bigl(u+2^{2/3}x\bigr)}{\mathrm {Ai}(u)^2}\,du,\quad x>0,
\end{equation}
and shift the integration path in the last integral to the path $P$ along the line, parallel to the imaginary axis, and running from $c-i\infty$ to $c+i\infty$, where $c=\tfrac132^{2/3}x$. A similar path was used in the proof of (ii) of Corollary 3.4 in \mycite{piet:89}. Asymptotically, as $x\to\infty$, the exponent is now given by
$$
\tfrac43(c+u)^{3/2}-\tfrac23(c+2^{2/3}x+u)^{3/2},
$$
which equals $-8x^{3/2}/\sqrt{27}$ at $u=0$. We get:
\begin{align*}
f_N(x)&=\frac1{2^{1/3}\pi i}\int_{c-i\infty}^{c+i\infty}\frac{\mathrm {Ai}\bigl(u+2^{2/3}x\bigr)}{\mathrm {Ai}(u)^2}\,du
\sim \frac{2c^{1/2}\sqrt{\pi}e^{-8x^{3/2}/\sqrt{27}}}{2^{1/3}(c+2^{2/3}x)^{1/4}\pi}\int_{-\infty}^{\infty}\exp\left\{-\frac{3u^2\sqrt{3}}{8\left(2^{1/3}\sqrt{x}\right))}\right\}\,du\\
&=\frac{2^{1/3}x^{1/4}e^{-8x^{3/2}/\sqrt{27}}}{3^{1/4}}\cdot\frac{2^{5/3}x^{1/4}}{3^{3/4}}=\frac{4e^{-8x^{3/2}/\sqrt{27}}\sqrt{x}}3, \quad x\to\infty.
\end{align*}
\eop

\bigskip
\begin{remark}\label{rem_meaning}
{\rm In the first proof, the integral in \eqref{Regtasymp} does not have a meaning for all complex values $t$. For example, we use values of $t$ with phase $\pi/3$ for which $t^3$ is negative. However, for all complex values of $t$ with phases in $[-\pi/3,\pi/3]$ we can give the integral in \eqref{Regtasymp} a meaning by turning the path of integration in the $w-$plane such that the integral remains convergent. For example Êwhen $\arg(t)=\pi/3$ we can take $\arg(w)=-\pi/2$.
}
\end{remark}

For two-sided Brownian motion we get similarly:

\begin{corollary}
\label{cor:2-sided}
Let $M$ be defined by
$$
M=\max_{t\in\R}\{W(t)-t^2\},
$$
where $W$ is standard two-sided Brownian motion, originating from zero.
and let $f_M$ and $F_M$ be the density and the distribution function of $M$, respectively. Then:
\begin{enumerate}
\item[(i)]
$$
f_M(x)=2f_N(x)F_N(x)\sim \frac{8\sqrt{x}}{3}\exp\left\{-\frac{8x^{3/2}}{3\sqrt{3}}\right\},\,x\to\infty.
$$
\item[(ii)]
$$
1-F_M(x)\sim \frac2{\sqrt{3}}\exp\left\{-\frac{8x^{3/2}}{3\sqrt{3}}\right\},\,x\to\infty.
$$
\end{enumerate}
\end{corollary}

\noindent
{\bf Proof.} This follows from Corollary 2.2 of \mycite{piet:10}, which gives the representation:
$$
f_M(x)=2f_N(x)F_N(x),\,x\ge0.
$$
\eop

\section{Concluding remarks}
\label{section:conclusion}
As pointed out to us by Svante Janson, the result implies certain facts for the moments of the distribution. For example, applying Theorem 4.5 in \mycite{janson:04} together with Theorem \ref{th:main_result} of the present paper gives:
$$
(E M^r)^{1/r}\sim \frac12 (3/2)^{1/3} (r/e)^{2/3}, r\to\infty.
$$

Theorem \ref{th:main_result} can also easily be extended to a result for
$$
M_c=\max_{t\ge0}\{W(t)-ct^2\},
$$
by using the scaling relation:
$$
M_c=c^{-1/3}M,
$$
see (1.7) of \mycite{janson:10}. So, for example, Theorem \ref{th:main_result} implies:
$$
\P\left\{M_c\ge x\right\}=\P\left\{M\ge c^{1/3}x\right\}\sim
\frac1{\sqrt{3}}\exp\left\{-\frac{8x^{3/2}\sqrt{c}}{3\sqrt{3}}\right\},\,x\to\infty.
$$
Also, by Lemma \ref{lemma:dens_representation}, the density of $M_c$ is given by:
\begin{equation}
\label{density_representation4}
f_{M_c}(x)=\frac{(4c)^{1/3}}{\pi}\,\mbox{\rm Re}\left(\int_0^{\infty}\frac{\ai\bigl(iu+(4c)^{1/3}x\bigr)}{\ai(iu)^2}\,du\right),\,x>0.
\end{equation}

\section{Appendix. Computing more coefficients of the asymptotic expansions}
\label{more}
Recall, see \eqref{defg2},
$$
g_2(t)=\frac{1}{\pi}\int_{\calP_2}\frac{e^{-itu}}{\ai(iu)^2}\,du.
$$
where we assume for the time being $t > 0$. We substitute $u=t^2 v$ and  write $\phi(u)=-t^3\psi(v)$, where $\phi(u)$ is given in \eqref{phidef}. This gives
$$
g_2(t)=\frac{t^2}{\pi}\int_{\widetilde{P}_2}\frac{e^{-t^3\psi(v)}}{e^{\tfrac43(it^2v)^{3/2}}\,\ai(it^2v)^2}\,dv,
$$
where the modified path $\widetilde{P}_2$ runs from $-i/4$ to $-i/4+\infty$, and
$$
\psi(v)= -\tfrac43(iv)^{3/2}+iv.
$$
We have $\psi(-i/4)=\frac{1}{12}$ and substitute
\begin{equation}
\label{trans}
\psi(v)-\tfrac1{12}=w^2.
\end{equation}
Locally, we prescribe for this mapping 
$$
w= (v+i/4)+\bigO\left((v+i/4)^2\right), \quad v\to-i/4.
$$
Upon inverting we have $v=\sum_{k=0}^\infty d_kw^k$, and the first coefficients are
$$
d_0=-\tfrac14i,\quad 
d_1=1, \quad
d_2=\tfrac{1}{3}i, \quad
d_3=\tfrac{2}{9},\quad
d_4=-\tfrac{8}{27}i,\quad
d_5=-\tfrac{14}{27}.
$$ 

The transformation \eqref{trans} gives
\begin{equation}
\label{trans2}
g_2(t)=t^2e^{-\frac1{12}t^3}\int_0^\infty e^{-t^3w^2} h(w)\,dw,
\end{equation}
where
$$
h(w)=\frac{1}{\pi e^{\frac43(it^2v)^{3/2}}\,\ai(it^2v)^2}\ \frac{dv}{dw}.
$$
Because the argument of the Airy function is large for large values of $t$ (for all values of $w\ge0$) we expand the function $h(w)$ first for large values of $t$. For this we need the well-known expansion  (see \mycite{olver:10}\footnote{http://dlmf.nist.gov/9.7.E5}):
\begin{equation}
\label{Airy_asymp1}
\ai(z)\sim\frac{e^{{-\zeta}}}{2\sqrt{\pi}z^{{1/4}}}
\sum _{{k=0}}^{{\infty}}(-1)^{k}\frac{u_{k}}{\zeta^{k}},
\quad z\to\infty, 
\quad \vert\phase\,z\vert<\pi,
\end{equation}
where 
$$
\zeta=\tfrac23z^{3/2},
$$
and
$$
u_{k}=\frac{(2k+1)(2k+3)(2k+5)\cdots(6k-1)}{(216)^{k}\,k!}=\frac{\Gamma(3k+\frac12)}{54^k\,k!\,\Gamma(k+\frac12)},
$$
The first coefficients are
$$
u_0=1,\quad 
u_1=\tfrac{5}{72}, \quad
u_2=\tfrac{385}{10368}, \quad
u_3=\tfrac{85085}{2239488}, \quad
u_4=\tfrac{37182145}{644972544}. 
$$ 

For $h(w)$ of \eqref{trans2} we obtain
$$
h(w)\sim 4t\sqrt{iv} \frac{dv}{dw}\sum _{{k=0}}^{{\infty}}\frac{w_{k}}{\zeta^{k}},\quad \zeta=\tfrac23t^3(iv)^{3/2},
$$
where the first coefficients are given by
$$
w_0=1,\quad 
w_1=\tfrac{5}{36}, \quad
w_2=-\tfrac{155}{2592}, \quad
w_3=\tfrac{17315}{279936}, \quad
w_4=-\tfrac{3924815}{40310784}. 
$$
We expand $h(w)=\sum_{k=0}^\infty h_k(t)w^k$ and observe that the terms with odd index are purely imaginary and those with even index are real. We obtain
\begin{equation}
\label{trans3}
\mbox{\rm Re}\left(g_2(t)\right)\sim t^2e^{-\frac1{12}t^3}\sum_{k=0}^\infty h_{2k}(t)\int_0^\infty e^{-t^3w^2}w^{2k}\,dw=
\tfrac12\sqrt{\pi}\,\sqrt{t}e^{-\frac1{12}t^3}\sum_{k=0}^\infty \frac{(\frac12)_k}{t^{3k}}h_{2k}(t).
\end{equation}

Each $h_{2k}(t)/t$ is an infinite series in negative powers of $t$, representing an asymptotic expansion of these coefficients. For example,
\begin{align*}
h_0(t)&\sim t\left(2+\tfrac{10}{3}t^{-3}-\tfrac{155}{9}t^{-6}+\tfrac{17315}{81}t^{-9}+\ldots,\right)\\
h_2(t)&\sim t\left( \tfrac43-\tfrac{340}{9}t^{-3}+\tfrac{27590}{27}t^{-6}-\tfrac{7445450}{243}t^{-9}+\ldots,\right)\\
h_4(t)&\sim t\left(-\tfrac{140}{27}+\tfrac{36020}{81}t^{-3}-\tfrac{6517750}{243}t^{-6}+\tfrac{3140421550}{2187}t^{-9}+\ldots,\right).
\end{align*}
Using these expansions and rearranging the series, we obtain
\begin{equation}
\label{trans4}
\mbox{\rm Re}\left(g_2(t)\right)\sim 
\sqrt{\pi}t^{3/2}e^{-\frac1{12}t^3}\sum_{k=0}^\infty \frac{a_k}{t^{3k}},
\end{equation}
where
\begin{align*}
\label{trans4}
&a_0=  1,\quad 
a_1=  2,\quad 
a_2=  -20,\quad 
a_3=  560,\quad
a_4=  -25520,\\
&a_5=  1601600,\quad 
a_6=  -127568000,\quad 
a_7=  12287436800.
\end{align*}

We have derived the expansion Êin \eqref{trans4} with the assumption $t>0$. 
As observed in Remark~\ref{rem_meaning} the integrals in \eqref{trans3} have a meaning for all phases of $t\in[-\pi/3,\pi/3]$. Therefore, we can use the result \eqref{trans4} in \eqref{asymp_representation} and obtain
\begin{equation}
\label{asymp_representation3}
f_N(x)\sim2^{2/3}\sqrt{\pi}\sum_{k=0}^\infty a_k \Phi_k(x),
\end{equation}
where
$$
\Phi_k(x)=\frac1{2\pi i}
\int_{\infty e^{-\pi i/3}}^{\infty e^{\pi i/3}} e^{\frac14t^3-2^{2/3}xt}t^{\frac32-3k}\,dt.
$$
We continue the analysis considering the Airy-type integrals
\begin{equation}
\label{aialpha}
A_\a(z) =\frac1{2\pi i}
\int_{\infty e^{-\pi i/3}}^{\infty e^{\pi i/3}} e^{\frac13w^3-zw}w^\a\,dw,
\end{equation}
for large values of $z$. The contour should cross the real axis for positive values of $w$. 
We have
$$
\Phi_k(x)=(4/3)^{5/6-k}A_{3/2-3k}\left(2^{4/3}3^{-1/3}x\right),\quad k=0,1,2,\ldots.
$$

For nonnegative integer values of $\a$ the functions $A_\a(z)$ are derivatives of the Airy function:
$$
A_n(z) =(-1)^n\frac{d^n}{dz^n}\ai(z),\quad n=0,1,2,\ldots.
$$
For negative integer  values of $\a$ they can be viewed as antiderivatives. We have, for example,
$$
A_{-1}(z) =\int_z^\infty \ai(t)\,dt.
$$

By using asymptotic expansions of the functions $A_\a(z)$ (as given in the lemma below) and the corresponding expansions for the $\Phi_k$ in \eqref{asymp_representation3} we finally obtain after a few manipulations of series
\begin{equation}
\label{fNxasymp}
f_N(x)\sim\frac{4\sqrt{x}}{3}\exp\left(-\frac{8x^{3/2}}{3\sqrt{3}}\right)\sum_{k=0}^\infty \frac{f_k}{x^{\frac32k}},\quad x\to\infty,
\end{equation}
where the first coefficients are
$$
f_0=  1,\quad 
f_1=  \tfrac{19}{48}\sqrt{3},
\quad 
f_2=  -\tfrac{3851}{1536},\quad 
f_3=   \tfrac{3380005}{221184}\sqrt{3},\quad
f_4=-\tfrac{6474441455}{14155776}.
$$

\begin{lemma}
\label{Aa}
The function  $A_\a(z)$ defined in \eqref{aialpha} has the following asymptotic expansion:
\begin{equation}
\label{Aalasymp}
A_\a(z)\sim\frac{z^{\a/2-1/4}e^{{-\zeta}}}{2\sqrt{\pi}}
\sum _{{k=0}}^{{\infty}}(-1)^{k}\frac{u_{k}^{(\a)}}{\zeta^{k}},
\quad z\to\infty, 
\quad \vert\phase\,z\vert<\pi, 
\end{equation}
where 
$$
\zeta=\tfrac23z^{3/2},
$$
and, in terms of the Gauss hypergeometric function,
\begin{equation}
\label{uk}
u_{k}^{(\a)}=\frac{\Gamma(3k+\frac12)}{54^{k}k!\Gamma(k+\frac12)}{}_2F_1\left(-\a,-2k;\tfrac12-3k;3\right).
\end{equation}
\end{lemma}

\noindent
{\bf Proof}. \ Substitute in  \eqref{aialpha} $w=\sqrt{z}s$ and write the integral in the form
$$
A_\a(z) =\frac{z^{\frac12(\a+1)}e^{-\zeta}}{2\pi i}
\int_{\infty e^{-\pi i/3}}^{\infty e^{\pi i/3}} e^{z^{3/2}\phi(s)}s^\a\,ds,
$$
where
$$
 \phi(s)=\tfrac13s^3-s+\tfrac23=\tfrac13(s+2)(s-1)^2.
$$
The substitution $t=(s-1)\sqrt{(s+2)/3)}$ gives
$$
A_\a(z) =\frac{z^{\frac12(\a+1)}e^{-\zeta}}{2\pi i}
\int_{-i\infty }^{i\infty } e^{z^{3/2}t^2}\psi(t)\,dt,\quad  \psi(t)=s^\a\frac{ds}{dt}.
$$
By expanding $\psi(t)=\sum_{k=0}^\infty \psi_k t^k$,
$$
A_\a(z) \sim\frac{z^{\a/2-1/4}e^{-\zeta}}{2\pi}
\sum_{k=0}^\infty (-1)^k \psi_{2k}\frac{\Gamma(k+\tfrac12)}{z^{3k/2}}.
$$
The coefficients $\psi_k$ can be represented as a Cauchy integral:
$$
\psi_k=\frac{1}{2\pi i}\int_C \psi(t) \frac{dt}{t^{k+1}},
$$
where $C$ is a small circle around the origin. By integrating in the $s-$plane, it follows that
$$
\psi_k=\frac{1}{2\pi i}\int_C s^\a \frac{ds}{((s+2)/3)^{(k+1)/2}(s-1)^{k+1}},
$$
where  $C$ is a circle around $s=1$ with radius less than 1. By using \mycite{olde:10}\footnote{http://dlmf.nist.gov/15.6.E3}:
$$
\psi_k=\frac{(-1)^k}{3^k k!}\frac{\Gamma(\frac32k+\frac12-\a)}{\Gamma(\frac12k+\tfrac12-\a)}\ 
{}_2F_1\left(-\a,-k;\tfrac12 k+\tfrac12-\a;-2\right).
$$
This form of the Gauss hypergeometric function has to be modified when $\tfrac12(k+1)=\a$, but using 
 \mycite{olde:10}\footnote{http://dlmf.nist.gov/15.8.E4}, and taking into account that $k$ is an integer, we can write
$$
{}_2F_1\left(-\a,-k;\tfrac12 k+\tfrac12-\a;-2\right)=
\frac{\Gamma(\frac32 k+\frac12)\,\Gamma(\frac12k+\frac12-\a)}
{\Gamma(\frac12 k+\frac12)\,\Gamma(\frac32 k+\frac12-\a)}\ 
{}_2F_1\left(-\a,-k;\tfrac12 -\tfrac32 k;3\right).
$$
This gives
$$
\psi_{2k}=\frac{1}{3^{2k} (2k)!}\frac{\Gamma(3k+\frac12)}{\Gamma(k+\tfrac12)}\ 
{}_2F_1\left(-\a,-2k; \tfrac12-3k;3\right).
$$
Using the duplication formula $(2k)!=2^{2k}k!\,\Gamma(k+\frac12)/\sqrt{\pi}$ we have proved the lemma.

\bigskip
\begin{remark}
{\rm The Gauss hypergeometric function in \eqref{uk} has a simple finite representation because $k$ is an integer.  Also, for $\a=1,2,3,\ldots$, giving results for the derivatives  of the Airy function, we have simple relations. For example, when $\a=1$ we have $u_k^{(1)}=u_k(1+6k)/(1-6k)=-v_k$, and the $v_k$ are used in the expansion of $\ai^\prime(z)$, see \mycite{olver:10}\footnote{http://dlmf.nist.gov/9.7.E6}.
}
\end{remark}

\begin{remark}
{\rm Similar expansion as in \eqref{Aalasymp} can be found in \mycite{drazin:81} (pp.~465--478) for slightly different functions.
}
\end{remark}

\vspace{0.3cm}
\noindent
{\bf Acknowledgements}
We thank the referees  for their careful reading and useful comments. We also want to thank Svante Janson and Jon Wellner for their useful comments and Zakhar Kabluchko for bringing  \mycite{huesler:99} to our attention. Nico M.~Temme acknowledges financial support from {\emph{Ministerio de Educaci\'on y Ciencia}}, project MTM2006--09050, {\emph{Ministerio de Ciencia e Innovaci\'on}}, project MTM2009-11686, and {\it Gobierno of Navarra}, Res. 07/05/2008.

\end{document}